\documentclass{article}

\usepackage{amsmath,amssymb,amsthm}
\usepackage{graphicx}
\usepackage{caption}
\usepackage{subcaption}
\usepackage{url}
\usepackage{a4wide}

\newtheorem{theorem}{Theorem}

\theoremstyle{definition}
\newtheorem{definition}[theorem]{Definition}

\begin{document}

\title{Combining the regularization strategy and the SQP to solve MPCC \\ - a MATLAB implementation\thanks{This is a preprint of a paper whose
final and definite form is in Journal of Computational
and Applied Mathematics.
Please cite this paper as: \emph{M. Teresa T. Monteiro and Helena Sofia Rodrigues (2011). Combining the regularization
strategy and the SQP to solve MPCC - a MATLAB implementation, Journal of Computational
and Applied Mathematics,
235(18), 5348--5356.}} }

\author{M. Teresa T. Monteiro $^{1,2}$\\
{\tt \small tm@dps.uminho.pt}
 \and Helena Sofia Rodrigues$^{1,3}$\\
{\tt \small sofiarodrigues@esce.ipvc.pt}
}


\date{$^1$ ALGORITMI Research Centre,University of Minho, Portugal\\[0.3cm]
$^2$ Department of Production and Systems, University of Minho, Portugal\\[0.3cm]
$^3$ Business School, Viana do Castelo Polytechnic Institute,\\
Portugal\\[0.3cm]
}

\maketitle



\begin{abstract}
Mathematical Program with Complementarity Constraints (MPCC) plays a very important role in many fields such as engineering design, economic equilibrium, multilevel game, and mathematical programming theory itself. In theory its constraints fail to satisfy a standard constraint qualification  such as the linear independence constraint qualification (LICQ) or the Mangasarian-Fromovitz constraint qualification (MFCQ) at any feasible point. As a result, the developed nonlinear programming theory may not be applied to MPCC class directly. Nowadays, a natural and popular approach is try to find some suitable approximations of an MPCC so that it can be solved by solving a sequence of nonlinear programs.

This work aims to solve the MPCC using nonlinear programming techniques, namely the SQP and the regularization scheme.  Some algorithms with two iterative processes, the inner and the external, were developed. A set of AMPL problems from MacMPEC  database \cite{MacMPEC} were tested. The algorithms performance comparative analysis was carried out.

\smallskip

\end{abstract}

%
%

\section{Introduction}
Mathematical Program with Complementarity Constraints is an exciting new application of nonlinear programming techniques working like a challenge for the scientific community. There exist several MPCC application areas like Engineering, Economics, Ecology among others. In Engineering one can distinguish the contact, obstacle and friction problems, process modeling, deformation and traffic congestion. MPCC problems arises in Economics related to game theory models like Nash and Stackelberg equilibrium, finances and taxes, markets competition. In Ecology the Quioto protocol is a real situation that can be modelled as a MPCC problem. Ralph \cite{Ralph} presents some MPCC applications like toll design in traffic networks or communication networks. The researchers have been spent lots of efforts studying  the MPCC theory and proposing different algorithms to solve MPCC efficiently. One can emphasize the work of Fukushima and Pang \cite{FukushimaPang}, Scholtes \cite{Scholtes01},  Anitescu \cite{Anitescu2005}, Scheel and Scholtes  \cite{ScheelScholtes}, Ralph and Wright \cite{RalphWright} and Fletcher \emph{et al.} \cite{FletcherLeyfferRalphScholtes}. The interior point method (IPM), the sequential quadratic programming (SQP), the smooth nonlinear programming, the penalty technique and regularization scheme are some strategies that have been studied to implement numerical algorithms. There also exist a growing collections of test problems.

An important reason why complementarity optimization problems are so pervasive in Engineering and Economics is because the concept of complementarity is synonymous with the notion of system equilibrium. This optimization problem is very difficult to solve because the usual constraint qualifications, necessary to guarantee the algorithms convergence, fail in all feasible points. This complexity is caused by the disjunctive constraints which lead to some challenging issues that typically are the main concern in the design of efficient solution algorithms. From the geometric point of view, its feasible region is not convex and not connected even in general.
Recently, it has been shown that MPCC can be solved efficiently and reliably as nonlinear program (NLP).  However this reformulation still continues to violate at any feasible point the same constraint qualifications (MFCQ), \emph{ie},  has no feasible point that satisfies the inequalities strictly.

 Recent studies of Scheel and Scholthes  \cite{ScheelScholtes} have proved that the strong stationarity of an MPCC is equivalent to the first order optimality conditions of the NLP equivalent. This fact motivated the cientifique community to use NLP approaches to deal with MPCC. Fletcher \emph{et al.} \cite{FletcherLeyfferRalphScholtes} complements these numerical observation giving a theoretical explanation for the good performance of the SQP method - they show that SQP is guaranteed to converge quadratically near a stationary point under relatively mild conditions.  Ralph and Wright \cite{RalphWright} described some properties of penalized and regularized nonlinear programming formulations of MPCC.  Based on these results we propose a general algorithm, and some alternatives, combining the SQP and the regularization strategy.

This paper is organized as follows. The next section defines the MPCC
problem. Some concepts related to the optimality conditions are
presented in Section 3. Section 4 introduces the regularization scheme, some alternative regularized formulations and some convergence issues. The implemented algorithms in MATLAB environment are detailed
in Section 5. Section 6 reports the MATLAB-AMPL interface, some numerical experiments to test the algorithms, a performance profiles analysis and some conclusions and future work ideas.

\section{Problem definition}
We consider Mathematical Program with Complementarity Constraints (MPCC) of the form
\begin{equation} \tag{MPCC}\label{mpcc}
\begin{tabular}{ll}
$min$ & $f(x)$ \\
$s.t.$ & $c_{i}(x)=0,\; i \in E,$ \\
& $c_{i}(x)\geq 0,\; i \in I,$\\
& $0\leq x_1\perp x_2\geq 0,$%
\end{tabular}
\end{equation}
\noindent where $f$ and $c$ are the nonlinear objective function and the constraint functions, respectively, assumed to be twice continuously differentiable. $E$ and $I$ are two\, disjoined\, finite\, index \,sets  with \,cardinality \,$p$ and $m$, \,respectively.  $x=(x_0,x_1,x_2)$ is a decomposition of the variables into $x_0 \in \mathbb{R}^n$ (control variables) and $(x_1,x_2) \in \mathbb{R}^{2\,q}$ (state variables). $0\leq x_1\perp x_2\geq 0: \mathbb{R}^{2q}\rightarrow\mathbb{R}^{q}$ are the $q$ complementarity constraints. The notation $x_{1}\perp x_{2}$ means that $x_{1j}x_{2j}=0$, for $j=1,\ldots,q$, \emph{ie}, the complementarity condition owns the disjunctive nature - $x_{1j}=0$ or $x_{2j}=0$, for $j=1,\ldots,q$.
This formulation doesn't exclude complementarity constraints like $0\leq G(x)\perp H(x)\geq 0$. With this kind of complementarity constraints, the problem can be reformulated, by introducing the slack variables $x_1$ and $x_2$. Grouping  all the equality constraints in $c_i(x)=0,$ the complementarity constraints have the form $0\leq x_1\perp x_2\geq 0$ and the problem presents the formulation (\ref{mpcc}). In this formulation all the properties like constraint qualifications or second order conditions are persevered. This formulation makes easy the properties theoretical study.

\section{Optimality conditions}
This section introduces some concepts related to stationarity and second order conditions. The optimality concepts follow the development of \cite{FletcherLeyfferRalphScholtes} and the corresponding proves can be consulted in this work.
One attractive way of solving (\ref{mpcc}) is to replace the complementarity constraints by a set of nonlinear inequalities, such as $x_{1j}\,x_{2j}\leq 0, j=1,\ldots,q$, and then solve the equivalent nonlinear program (NLP):

\begin{equation}\label{nlp_eq}
\begin{tabular}{ll}
$min$ & $f(x)$ \\
$s.t.$ & $c_{i}(x)=0,\; i \in E,$ \\
& $c_{i}(x)\geq 0,\; i \in I,$\\
& $x_1\geq 0,\,x_2\geq 0,$\\
& $x_{1j}x_{2j}\leq 0,\, j=1,\ldots,q$.%
\end{tabular}
\end{equation}
It has been shown \cite{ScheelScholtes} that (\ref{nlp_eq}) violates the Mangasarian-Fromovitz constraint qualification (MFQC) at any feasible point. This failure of MFQC implies that multiplier set is unbounded, the central path fails to exist, the active constraints normals are linearly dependent, and the NLP linearizations can became inconsistent arbitrarily close to a solution. Recently, new developments motivated the interest in the analysis of NLP solvers applied to (\ref{nlp_eq}) based on the success of SQP methods - the simple observation of Scholtes is that strong stationarity is equivalent to the KKT conditions of (\ref{nlp_eq}). This fact implies the existence of bounded multipliers.

Consider  two index sets: $X_{1},X_{2}\subset
\{1,\ldots,q\}$ with $X_{1} \cup X_{2} = \{1,\ldots,q\}$, denoting the corresponding complements in $\{1,\ldots,q\}$ by  $X_{1}^{\perp}$ e $X_{2}^{\perp}$. For each pair of index one define the relaxed NLP corresponding to  (\ref{mpcc}):

\begin{equation} \tag{NLP-rel}\label{NLP-rel}
\begin{tabular}{ll}
$min$ & $f(x)$ \\
$s.t.$ & $c_i(x)=0,\; i \in E,$ \\
& $c_i(x)\geq 0, \; i \in I,$\\
& $x_{1j}= 0, \,\,\forall j \in X_{2}^{\perp},$\\
& $x_{2j}= 0, \,\,\forall j \in X_{1}^{\perp},$\\
& $x_{1j}\geq 0, \,\,\forall j \in X_{2},$\\
& $x_{2j}\geq 0, \,\,\forall j \in X_{1}.$\\
\end{tabular}
\end{equation}

Concepts like constraints qualification, stationarity and second order conditions of the MPCC problem will be defined in terms of (\ref{NLP-rel}).
The linear independence constraint qualification, LICQ, is extended to MPCC that is MPCC-LICQ:

\begin{definition}{MPCC-LICQ}

Consider $x_1,x_2 \geq 0$ and
define:
\[X_{1}=\{j : x_{1j}=0\},\]
\[X_{2}=\{j : x_{2j}=0\}.\]
The MPCC problem verifies  the MPCC-LICQ at $x$ if the corresponding (\ref{NLP-rel}) verifies the LICQ.
\end{definition}
If $x^{*}$ is a local solution of (\ref{NLP-rel}) and satisfies ${x_1^*}^T x_2^*=0$, then $x^{*}$ is also a local solution of original MPCC.

\noindent There are several kinds of stationarity defined for MPCC problem. Among them, the strong stationarity is the following one:

\begin{definition}{Strong stationarity} \label{forte}

$x^{*}$ is a strong stationary point if exist Lagrange multipliers $\lambda$, $\widehat{\nu}_1$ and
$\widehat{\nu}_2$ so that:
\begin{equation}\label{ooo}
\begin{tabular}{l}
$\nabla f^{*}-\left[\nabla (c_{i}^{*}), {\small i\in E}\;\;:\;\;\nabla (c_i ^{*}), {\small i \in I}\right]\lambda-\left(\begin{array}{c} 0 \\
\widehat{\nu}_1 \\
\widehat{\nu}_2\\
\end {array}\right)=0,$ \\
$c_i^*=0, i \in E,$ \\
$c_i^*\geq 0, i \in I, $\\
$x_{1}^{*}\geq 0, $\\
$x_{2}^{*} \geq 0, $\\
$x_{1j}^{*}= 0$ or $ x_{2j}^{*} = 0,$\\
$\lambda_{i}\geq 0, i \in I, $\\
$c_{i}\lambda_{i}= 0, $\\
$x_{1j}^{*} \widehat{\nu}_{1j}= 0, $\\
$x_{2j}^{*} \widehat{\nu}_{2j}= 0, $\\
if $x_{1j}^{*} = x_{2j}^{*}=0$ then $\widehat{\nu}_{1j}\geq 0$ and $\widehat{\nu}_{2j} \geq 0. $\\
\end{tabular}
\end{equation}
\end{definition}

Note that (\ref{ooo}) are the first order optimality conditions of the (\ref{NLP-rel}) at $x^{*}$.

At $x^{*}$ consider
\[A=\left[\nabla (c_i^{*}), {\small i\in E}\;\;:\;\;\nabla (c_i ^{*}), {\small i \in I}\cap A^*\;\;:\;\;\begin{array}{c} 0 \\
I_1^* \\
0\\
\end {array}\;\;:\;\;\begin{array}{c} 0 \\
0 \\
I_2^* \\
\end {array}\right]=:\left[a_i^*\right]_{i\in A^*},\]

\noindent where $I_1^*:=[e_i]_{i\in X_1^*}\;\; \text{ and } \;\;I_2^*:=[e_i]_{i\in X_2^*}$
are part of the $q \times q$ identity matrix  corresponding to the active simple bound constraints.
The set of feasible directions with null curvature of (\ref{NLP-rel}) is defined by:
\[
S^*=\{s| s \neq 0, \nabla (f^*)^{T} s=0,\; (a_{i}^*)^{T} s=0,\; i \in A_{+}^*,\; (a_{i}^*)^{T} s \geq 0,\; i\in A^{*}\backslash A_+ ^{*} \}.
\]
\noindent where $A^*$ is the index set of active constraints and  $A_+^* \subset A^*$ is the index set of  nondegenerated active constraints.

The second-order sufficient condition (SOSC) for MPCC is given as follows:
\begin{definition} {MPCC-SOSC}

A strong stationary point $x^{*}$ with multipliers $(\lambda^{*},\widehat{\nu}_1^{*},\widehat{\nu}_2^{*})$ verifies the a
MPCC-SOSC if all direction $s \in S^{*}$ satisfies
$s^{T}\nabla^{2}L^{*}s>0$ where $\nabla^{2}L^{*}$ represents the Hessian matrix of the Lagrangean function of (\ref{NLP-rel}) at $
(x^{*},\lambda^{*},\widehat{\nu}_1^{*},\widehat{\nu}_2^{*})$.
\end{definition}

\section{Regularization scheme}
The complementarity constraints are responsible for the main difficulties of an MPCC. In order to overcame this hard problem, some parameters are introduced to smooth or relax these constraints.
Ralph and Wright \cite{RalphWright} present several regularization schemes and the corresponding properties in order to solve MPCC. The same authors study a regularization scheme that is analyzed  by Scholtes \cite{Scholtes01}  where  (\ref{mpcc}) is approximated by the following  NLP problem with a non negative scalar parameter $t$ decreasing to zero:
\begin{equation}\label{Reg}
\begin{tabular}{lll}
$Reg(t)$: & $min$ & $f(x)$ \\
& $s.t.$ & $c_{i}(x)=0, i\in E,$ \\
& &$c_{i}(x)\geq 0, i\in I, $\\
& &$x_1\geq 0,\, x_2\geq 0,$\\
& &$x_{1j}x_{2j}\leq t,\,\,\; j=1,...,q.$
\end{tabular}
\end{equation}
The solution of this problem is denoted by $x(t)$. The idea is to find a local minimum $x_k$ of $Reg(t_k)$ where $0<t_k\rightarrow 0.$ Suppose $x^*$ is a limit point of $\{x_k\}$, then $x^*$ is feasible for $Reg(0)$ hence for the (MPCC). As the $Reg(0)$ is equivalent to (\ref{nlp_eq}), the regularization scheme can be used by applying a NLP algorithm to $Reg(t)$ for a sequence of problems where $t$ is positive  and tends to zero. In this sequence, the result of each minimization, \emph{ie}, the approximate minimizer of the original problem, is the initial approximation of the next minimization. This minimization sequence represents the external iterative process.

\subsection{Other alternative regularized formulations}
Beyond this previous formulation, other equivalent formulations of (MPCC) are implemented in this work. The first one replaces the complementarity constraints by only one constraint, for this reason it is named $Reg$-$one$:
\begin{equation} \label{Reg-one}
\begin{tabular}{lll}
$Reg$-$one(t)$: & $min$ & $f(x)$ \\
& $s.t.$ & $c_{i}(x)=0, i\in E,$ \\
& & $c_{i}(x)\geq0, i\in I, $\\
& & $x_{1}\geq 0,\,\,\,\,x_{2}\geq 0,$\\
& & $x_{1}^T x_{2}\leq t,$
\end{tabular}
\end{equation}
where $x_{1}^T x_{2}\leq t \Leftrightarrow \displaystyle {\sum_{k=1}^q } x_{1k}x_{2k} \leq t.$ This formulation is of interest in computation because it has fewer constraints than (MPCC). This formulation was studied by Scholtes \cite{Scholtes01} and Anitescu \cite{Anitescu2005}  and we refer the reader to these sources for a detailed theoretical analysis.

Another plausible regularization was proposed by Ralph \cite{RalphWright} where the inequalities of (\ref{Reg}) are replaced by equalities:
\begin{equation} \label{Reg-eq}
\begin{tabular}{lll}
$Reg$-$eq(t)$: & $min$ & $f(x)$ \\
& $s.t.$ & $c_{i}(x)=0, i\in E,$ \\
& &$c_{i}(x)\geq 0, i\in I, $\\
& &$x_1\geq 0,\, x_2\geq 0,$\\
& &$x_{1j}x_{2j}= t,\,\,\; j=1,...,q.$
\end{tabular}
\end{equation}

Based on the previous formulations, another regularization scheme is proposed in this work:
\begin{equation}\label{Reg-eq-one}
\begin{tabular}{lll}
$Reg$-$eq$-$one(t)$: & $min$ & $f(x)$ \\
& $s.t.$ & $c_{i}(x)=0, i\in E,$ \\
& & $c_{i}(x)\geq0, i\in I, $\\
& & $x_{1}\geq 0,\,\,\,\,x_{2}\geq 0,$\\
& & $x_{1}^T x_{2}= t,$
\end{tabular}
\end{equation}
in which the complementarity constraints have been replaced by only one equality constraint $x_{1}^T x_{2}=t \Leftrightarrow \displaystyle{\sum_{k=1}^q }x_{1k}x_{2k}=t.$

\section{MATLAB Algorithms}

\renewcommand{\baselinestretch}{1.3}
\textbf{General algorithm}   \\
\begin{tabular}{||l||}
\hline \hline
Initialization: $t_0$, $k=0$;\\
Tolerances: $t_{min},\, k_{max},\, \epsilon_1,\, \epsilon_2$;\\
Inner iterations counter: $it\_int=0$;\\
Problem information (amplfunc): $x_0,\, lb,\,ub,\, cl,\, cu,\, cv$;\\
Problem dimension: $n,\,m,\,p,\,q$;\\
\textbf{\textsc{REPEAT}}\\
\hspace{0.4cm}\textbf{Step 1} - Built the constraints;\\
\hspace{0.4cm}\textbf{Step 2} - Run the MATLAB function:\\
\hspace{0.4cm}$[x,f,LAMBDA, output]$=\texttt{fmincon('function',....'constraint')};\\
\hspace{0.4cm}\textbf{Step 3} - $it\_int \leftarrow it\_int + output.iterations$;\\
\hspace{0.4cm}\textbf{Step 4} - Approximation update: $x_{k+1}\leftarrow x$;\\
\hspace{0.4cm}\textbf{Step 5} - Lagrange multipliers update; \\
\hspace{0.4cm}\textbf{Step 6} - Update $t\, (0<\rho_2<1)$: $t_{k+1}\leftarrow t_{k}\times \rho_2$;\\
\hspace{0.4cm}\textbf{Step 7} - $k\leftarrow k+ 1$;\\
\textbf{\textsc{UNTIL}} Stop criterium\\
\hline \hline
\end{tabular}

 Four algorithms were implemented, $Reg$, $Reg$-$one$, $Reg$-$eq$ and $Reg$-$eq$-$one$, corresponding to (\ref{Reg}), (\ref{Reg-one}), (\ref{Reg-eq}) and (\ref{Reg-eq-one}) regularized formulations, respectively. The diference between them is the complementarity constraints treatment. For this reason only the general algorithm is reported. The algorithm has two iterative processes - the external one, performs a sequence of minimization problems. Each external iteration executes the  \texttt{fmincon} MATLAB subroutine, which implements the SQP strategy. The \texttt{fmincon} call is the inner iterative process.

Next, the meaning of some parameters and procedures is reported.
The initial value of the regularization parameter is $t_0$; $it\_int$ and $k$ are the inner and external iteration counter, respectively; $t_{min}$ and $k_{max}$ are the regularization parameter limit and the external iterations limit, respectively; $\epsilon_1$ and $\epsilon_2$ are small positive constants. The specific information about the test problem is  $x_0,\, lb,\,ub,\, cl,\, cu,\, cv$ and will be explained in the next section. Using this information the final dimensions are calculated - $n,\,m,\,p,\,q$.

The external iterative procedure starts at the Step 1 with the constraints treatment routine  - each algorithm performs its own strategy replacing the complementarity constraints. Step 2 refers to the inner iterative procedure - the \texttt{fmincon} call performs the SQP and its input parameters will be explained in the next section. Step 3 updates the inner iterations using an output \texttt{fmincon} parameter. The iteration $k$ of the external iterative procedure updates the approximation of the solution $x_{k+1}$ with the solution of the  inner iterative procedure.
In Step 5 the Lagrange multipliers of the original problem, MPCC,  are updated using the structure LAMBDA which is a \texttt{fmincon} output parameter. Step 6 updates the regularization parameter with a smaller value and Step 7 increments the external iteration counter. This external iterative procedure is controlled by a stop criterium that is the disjunction of the four conditions:
\begin{equation}\label{cp}
t\leq t_{min}\;\; \vee\;\;
k=k_{max}\;\;\vee\;\; \dfrac{\|x_{k+1}-x_k\|}{\|x_{k+1}\|}\leq
\epsilon_1 \;\; \vee\;\;\;\|\nabla \mathcal{L}(x,\lambda)\| \leq \epsilon_2
\end{equation}
where $\dfrac{\|x^{k+1}-x^k\|}{\|x^{k+1}\|}$ is the relative error estimation in $x$ and  $\|\nabla \mathcal{L}(x,\lambda)\|$ is he stationarity evaluation.

\subsection{Convergence issues}
In this work, the proposed algorithm is based on the ideas presented by Scholtes  \cite{Scholtes01} and Ralph and Wright \cite{RalphWright}. Scholtes (\cite{Scholtes01}, Theorem 4.1) shows that in the neighborhood of the solution $x^*$ of (\ref{mpcc}), satisfying certain assumptions, exists only a stationary point $x(t)$ for $Reg(t)$ for all $t$ positive value  sufficiently small and  furthermore verifies $||x(t)-x^*||=O(t)$. The convergence behavior of a sequence of stationary points of a parametric NLP which regularizes an MPCC in the form of complementarity conditions is presented. Some important convergence properties are also proved: accumulation points are feasible points of the MPCC; they are M-stationary if, in addition, an approaching subsequence satisfies second order necessary conditions, and they are B-stationary if, in addition, an upper level strict complementarity condition holds. The same article shows that every local minimizer of the MPCC which satisfies the linear independence, upper level strict complementarity, and a second order optimality condition can be embedded into a locally unique piecewise smooth curve of local minimizers of the parametric NLP.

Based on this work, Ralph and Wright \cite{RalphWright} show that $Reg(t)$ has a local solution, possibly not only, such that $O(t^{\frac{1}{2}})$. They also prove that the Lagrange multipliers of the $Reg(t)$ solution are bounded and satisfies $O(t)$. The authors describe some properties of the solutions to the regularized formulations $Reg(t)$ (\ref{Reg}) to the MPCC (\ref{mpcc}): distance between solutions of (\ref{Reg}) and (\ref{mpcc}), boundedness of Lagrange multipliers, local uniqueness and smoothness of the solution mapping, under some assumptions on (\ref{mpcc}) at a local solution $x^*$.

\section{Computational experiments}

This section summarizes the numerical experiences using AMPL test problems from MacMPEC \cite{MacMPEC}. Details on problem size and characteristics can be found there. The computational experiences were made on a \emph{centrino} with 504MB of RAM, Windows operating system and 95 problems  are used to test the four algorithms.

\subsection{AMPL-MATLAB interface}

The proposed algorithms were implemented in MATLAB language using the  \texttt{fmincon} routine from the MATLAB Optimization toolbox (version 7.0.1.). This routine finds a constrained minimum of a several variables function starting at an initial estimate using a SQP method.

The test problems \cite{MacMPEC} are in AMPL language \cite{AMPL} whose files have the \emph{.mod } extension. As the MATLAB only recognizes the data in the \emph{.nl} shape the files have to be converted from \emph{.mod } to \emph{.nl} format.
An interface between AMPL and MATLAB was developed. The algorithms use the objective and constraints derivatives provided by AMPL.

\begin{figure}[h]\label{interf2}
\centering
\includegraphics[scale=0.50]{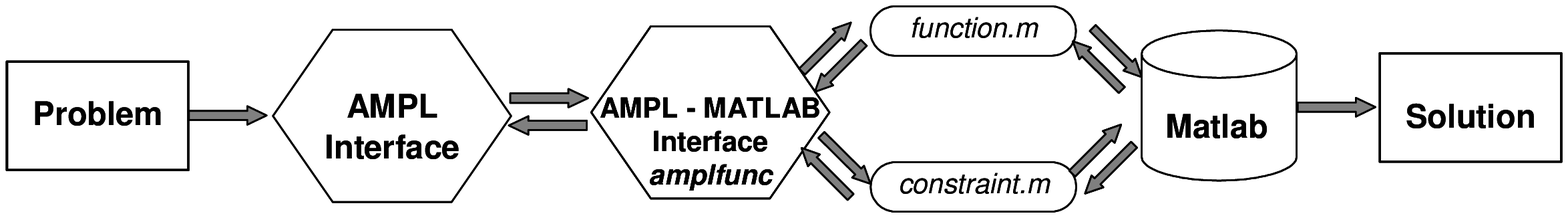}
\caption{AMPL-MATLAB interface}
\end{figure}

The MEX amplfunc function allows the MATLAB accessing to \emph{.mod} data. In order to make uniform the amplfunc function parameters two M-files (function.m and constraint.m) were developed to prepare the \texttt{fmincon} input parameters related to the objective and constraints functions. The  AMPL-MATLAB connection scheme is in Figure  \ref{interf2}.

To test a problem it is necessary the following input in the MATLAB command window:
\begin{verbatim}
>>[x0,lb,ub,v,cl,cu]=amplfunc('problem')
\end{verbatim}

The output information is: $x0$ - initial estimate, $lb$ and $ub$ - lower and upper limit of $x$ variable, $v$ dual problem initial estimate, $cl$ and $cu$ -  constraints lower and upper limit, $cv$ -  complementarity identifier array. The $cv$ array allows to identify the complementarity constraints:
\begin{equation*}
cv(i)=\left\{\begin{array}{ll}>0&\text{if the } i \text{ constraint}  \text{ complements with } x(cv(i))\\
0 & \text{ otherwise}.\end{array}\right.
\end{equation*}

\subsection{Numerical tests}

 The algorithms  were tested using two different update values for the regularization parameter: $t=0.1\times t$ performed better than $t=0.05 \times t$. The Table \ref{falhas} reports the results of the four algorithms with respected to their robustness.

\begin{table}
\begin{tabular}{|c||c|c|c|c|}
  \hline

  Algorithm & $Reg$ & $Reg$-$one$ & $Reg$-$eq$ & $Reg$-$eq$-$one$ \\  \hline \hline

   failures & 3 & 1 & 54 & 29 \\\hline
   \% & 3.2& 1.1 & 56.8 & 30.5  \\
  \hline
\end{tabular} \caption{\label{falhas}Failures}
\end{table}

The Algorithms $Reg$-$eq$ and  $Reg$-$eq$-$one$ present a significant number of failures and for this reason they are not considered in the following detailed results analysis.

The Tables \ref{tab1} and \ref{tab2} show the problem name, $n$ is the number of variables, $m$, $p$ and $q$ are the number of inequality, equality
and complementarity constraints, respectively. The next three columns report the Algorithm $Reg$ results - the minimum found ($f^*$), $it\_int$ and $it\_ext$ are the iteration counter for the inner and external iterations, respectively. The last three columns concern to the Algorithm $Reg$-$one$ and have the same meaning that the previous one. When the algorithm doesn't converge we denote by NC. To confirm the robustness of algorithms  a comparison with two codes (MacMPEC and Biegler   \cite{wachbieg304}) was performed.

\begin{table}[h]
\centering \caption{\label{tab1}Algorithm $Reg$ and Algorithm $Reg$-$one$}
\renewcommand{\baselinestretch}{1}
{\tiny
\begin{tabular}[c]{|l|c|c|c|c||r|c|c||r|c|c|}
\hline \noalign{\smallskip}
Name	&	n	&	m	&	p	&	q	&	f$^*$ ($Reg$)	&	it\_int	&	it\_ext	&	f$^*$ ($Reg$-$one$)	&	 it\_int	 &	it\_ext	\\
\noalign{\smallskip}\hline\noalign{\smallskip}
 bard1         	&	5	&	3	&	1	&	3	&	1.700000E+01	&	19	&	6	&	1.700000E+01	&	22	&	6	 \\
 bard1m        	&	6	&	3	&	1	&	3	&	1.700000E+01	&	20	&	7	&	1.700000E+01	&	23	&	7	 \\
 bard2         	&	12	&	4	&	5	&	3	&	6.163000E+03	&	4	&	2	&	6.163000E+03	&	4	&	2	 \\
 bard2m        	&	12	&	4	&	5	&	3	&	-6.598000E+03	&	6	&	2	&	-6.598000E+03	&	6	&	2	 \\
 bard3         	&	6	&	2	&	3	&	1	&	-1.267871E+01	&	7	&	6	&	-1,y267871E+01	&	7	&	6	 \\
 bard3m        	&	6	&	4	&	1	&	3	&	-1.267871E+01	&	28	&	7	&	-1.267871E+01	&	26	&	8	 \\
 bar-truss-3   	&	35	&	6	&	28	&	6	&	1.016657E+04	&	16	&	4	&	1.016657E+04	&	21	&	5	 \\
 bilevel1      	&	10	&	7	&	2	&	6	&	4.999999E+00	&	21	&	5	&	-2.812000E-07	&	25	&	6	 \\
 bilevel2      	&	16	&	9	&	4	&	8	&	-6.600000E+03	&	39	&	8	&	-6.600000E+03	&	21	&	8	 \\
 bilevel3      	&	11	&	4	&	6	&	3	&	-1.267871E+01	&	37	&	8	&	-1.267871E+01	&	42	&	8	 \\
 bilin         	&	8	&	7	&	0	&	6	&	-1.600000E-09	&	28	&	8	&	-9.400000E-09	&	18	&	8	 \\
 dempe         	&	3	&	1	&	1	&	1	&	3.125000E+01	&	60	&	8	&	3.125000E+01	&	60	&	8	 \\
 design-cent-1 	&	12	&	5	&	6	&	3	&	3.283210E-05	&	28	&	2	&	3.283210E-05	&	28	&	2	 \\
 design-cent-2 	&	13	&	9	&	6	&	3	&	NC	&		&		&	NC	&		&		\\
 design-cent-3 	&	15	&	5	&	6	&	3	&	3.208987E-01	&	24	&	3	&	1.691148E-01	&	53	&	2	 \\
 design-cent-4 	&	22	&	10	&	10	&	8	&	0.000000E+00	&	4	&	3	&	0.000000E+00	&	4	&	3	 \\
 desilva       	&	6	&	2	&	2	&	2	&	-1.000000E+00	&	5	&	4	&	-1.000000E+00	&	5	&	4	 \\
 df1           	&	2	&	3	&	0	&	1	&	0.000000E+00	&	2	&	2	&	0.000000E+00	&	2	&	2	 \\
 ex9.1.1       	&	13	&	5	&	7	&	5	&	-1.300000E+01	&	10	&	5	&	-1.300000E+01	&	11	&	6	 \\
 ex9.1.2       	&	8	&	2	&	5	&	2	&	-6.250000E+00	&	8	&	6	&	-6.250000E+00	&	9	&	7	 \\
 ex9.1.3       	&	23	&	6	&	15	&	6	&	-2.920000E+01	&	41	&	7	&	-2.920000E+01	&	23	&	7	 \\
 ex9.1.4       	&	8	&	2	&	5	&	2	&	-3.700000E+01	&	8	&	5	&	-3.700000E+01	&	8	&	5	 \\
 ex9.1.5       	&	13	&	5	&	7	&	5	&	-1.000000E+00	&	16	&	7	&	-1.000000E+00	&	15	&	8	 \\
 ex9.1.6       	&	14	&	6	&	7	&	6	&	-4.900000E+01	&	38	&	6	&	-4.900000E+01	&	31	&	6	 \\
 ex9.2.1       	&	10	&	4	&	5	&	4	&	1.700000E+01	&	40	&	6	&	1.700000E+01	&	81	&	7	 \\
 ex9.2.2       	&	9	&	4	&	4	&	3	&	9.999984E+01	&	60	&	8	&	9.999972E+01	&	63	&	8	 \\
 ex9.2.3       	&	14	&	5	&	8	&	4	&	5.000000E+00	&	5	&	2	&	5.000000E+00	&	5	&	2	 \\
 ex9.2.4       	&	8	&	2	&	5	&	2	&	5.000000E-01	&	37	&	8	&	5.000000E-01	&	29	&	8	 \\
 ex9.2.5       	&	8	&	3	&	4	&	3	&	9.000000E+00	&	37	&	8	&	9.000000E+00	&	48	&	8	 \\
 ex9.2.6       	&	16	&	6	&	6	&	6	&	-1.000000E+00	&	20	&	8	&	-1.000000E+00	&	20	&	8	 \\
 ex9.2.7       	&	10	&	4	&	5	&	4	&	1.700000E+01	&	40	&	6	&	1.700000E+01	&	81	&	7	 \\
 ex9.2.8       	&	6	&	2	&	3	&	2	&	1.500000E+00	&	15	&	7	&	1.500000E+00	&	14	&	7	 \\
 ex9.2.9       	&	9	&	3	&	5	&	3	&	2.000000E+00	&	8	&	3	&	2.000000E+00	&	8	&	3	 \\
 ex9-1-7n      	&	17	&	6	&	9	&	6	&	-2.300000E+01	&	35	&	7	&	-2.600000E+01	&	23	&	7	 \\
 ex9-1-9n      	&	12	&	5	&	6	&	5	&	3.111111E+00	&	24	&	7	&	3.111111E+00	&	25	&	7	 \\
 ex9-1-10n     	&	11	&	4	&	5	&	3	&	-3.250000E+00	&	9	&	6	&	-3.250000E+00	&	10	&	7	 \\
 gauvin        	&	3	&	2	&	0	&	2	&	2.000000E+01	&	22	&	6	&	2.000000E+01	&	20	&	6	 \\
 gnash10       	&	13	&	8	&	4	&	8	&	-2.308232E+02	&	25	&	7	&	-2.308232E+02	&	21	&	5	 \\
 gnash11       	&	13	&	8	&	4	&	8	&	-1.299119E+02	&	23	&	6	&	-1.299119E+02	&	21	&	5	 \\
 gnash12       	&	13	&	8	&	4	&	8	&	-3.693311E+01	&	19	&	6	&	-3.693311E+01	&	21	&	6	 \\
 gnash13       	&	13	&	8	&	4	&	8	&	-7.061784E+00	&	20	&	5	&	-7.061783E+00	&	23	&	8	 \\
 gnash14       	&	13	&	8	&	4	&	8	&	-1.790463E-01	&	30	&	8	&	-1.790463E-01	&	26	&	8	 \\
 gnash15       	&	13	&	8	&	4	&	8	&	-3.546991E+02	&	117	&	6	&	-3.546991E+02	&	33	&	6	 \\
 gnash16       	&	13	&	8	&	4	&	8	&	-2.414420E+02	&	117	&	6	&	-2.414420E+02	&	25	&	6	 \\
 gnash17       	&	13	&	8	&	4	&	8	&	NC	&		&		&	-9.074910E+01	&	36	&	8	\\
 gnash18       	&	13	&	8	&	4	&	8	&	NC	&		&		&	-2.569822E+01	&	56	&	8	\\
 gnash19       	&	13	&	8	&	4	&	8	&	-6.116708E+00	&	43	&	7	&	-6.116708E+00	&	51	&	8	 \\
 hs044-i       	&	20	&	10	&	4	&	10	&	1.561777E+01	&	34	&	8	&	2.061065E+01	&	29	&	8	 \\
 incid-set1-8  	&	117	&	70	&	49	&	49	&	0.000000E+00	&	7	&	2	&	0.000000E+00	&	7	&	2	 \\
 incid-set1c-8 	&	117	&	77	&	49	&	49	&	0.000000E+00	&	7	&	2	&	0.000000E+00	&	7	&	2	 \\
 incid-set2-8  	&	117	&	70	&	49	&	49	&	5.075396E-03	&	99	&	8	&	5.095247E-03	&	44	&	8	 \\
 incid-set2c-8 	&	117	&	77	&	49	&	49	&	5.642375E-03	&	59	&	8	&	5.640741E-03	&	106	&	8	 \\
 jr1           	&	2	&	1	&	0	&	1	&	5.000000E-01	&	3	&	2	&	5.000000E-01	&	3	&	2	 \\
 jr2           	&	2	&	1	&	0	&	1	&	5.000000E-01	&	21	&	8	&	5.000000E-01	&	21	&	8	 \\
 kth1          	&	2	&	1	&	0	&	1	&	0.000000E+00	&	3	&	2	&	0.000000E+00	&	3	&	2	 \\
 kth2          	&	2	&	1	&	0	&	1	&	0.000000E+00	&	3	&	2	&	0.000000E+00	&	3	&	2	 \\
 kth3          	&	2	&	1	&	0	&	1	&	0.000000E+00	&	1	&	1	&	0.000000E+00	&	1	&	1	 \\
 liswet1-050   	&	152	&	51	&	52	&	50	&	1.399428E-02	&	4	&	2	&	1.399428E-02	&	4	&	2	 \\
 nash1         	&	6	&	2	&	2	&	2	&	2.563700E-06	&	30	&	8	&	1.138420E-05	&	25	&	8	 \\
 \noalign{\smallskip} \hline
\end{tabular}
}
\end{table}

\begin{table}[h]
\centering \caption{\label{tab2} Algorithm $Reg$ and Algorithm $Reg$-$one$ (cont.)}
\renewcommand{\baselinestretch}{1}
{\tiny
\begin{tabular}[c]{|l|c|c|c|c||r|c|c||r|c|c|}
\hline \noalign{\smallskip}
\hline \noalign{\smallskip}
Name	&	n	&	m	&	p	&	q	&	f$^*$ ($Reg$)	&	it\_int	&	it\_ext	&	f$^*$ ($Reg$-$one$)	&	 it\_int	 &	it\_ext	\\
\noalign{\smallskip}\hline\noalign{\smallskip}

 outrata31     	&	5	&	4	&	0	&	4	&	3.207700E+00	&	23	&	7	&	3.207700E+00	&	23	&	7	 \\
 outrata32     	&	5	&	4	&	0	&	4	&	3.449404E+00	&	29	&	8	&	3.449404E+00	&	30	&	7	 \\
 outrata33     	&	5	&	4	&	0	&	4	&	4.604254E+00	&	25	&	7	&	4.604254E+00	&	29	&	8	 \\
 outrata34     	&	5	&	4	&	0	&	4	&	6.592684E+00	&	29	&	8	&	6.592684E+00	&	29	&	7	 \\
 pack-comp1-8  	&	107	&	72	&	49	&	49	&	6.000000E-01	&	19	&	8	&	6.000000E-01	&	18	&	8	 \\
 pack-comp1c-8 	&	107	&	79	&	49	&	49	&	6.000000E-01	&	15	&	8	&	6.000000E-01	&	17	&	8	 \\
 pack-comp2-8  	&	107	&	72	&	49	&	49	&	6.731171E-01	&	21	&	8	&	6.731171E-01	&	22	&	8	 \\
 pack-comp2c-8 	&	107	&	79	&	49	&	49	&	6.734582E-01	&	9	&	2	&	6.734582E-01	&	9	&	2	 \\
 pack-rig1-8   	&	87	&	40	&	46	&	32	&	7.879311E-01	&	24	&	8	&	7.879318E-01	&	25	&	8	 \\
 pack-rig1c-8  	&	87	&	47	&	46	&	32	&	7.882998E-01	&	20	&	8	&	7.883001E-01	&	23	&	8	 \\
 pack-rig1p-8  	&	105	&	55	&	49	&	47	&	7.879301E-01	&	24	&	8	&	7.879318E-01	&	28	&	8	 \\
 pack-rig2-8   	&	85	&	38	&	46	&	30	&	7.804042E-01	&	23	&	8	&	7.804027E-01	&	25	&	8	 \\
 pack-rig2c-8  	&	85	&	45	&	46	&	30	&	7.993058E-01	&	21	&	8	&	7.993051E-01	&	21	&	8	 \\
 pack-rig2p-8  	&	103	&	53	&	49	&	45	&	7.804042E-01	&	35	&	8	&	7.804027E-01	&	39	&	8	 \\
 portfl-i-1    	&	87	&	1	&	13	&	0	&	1.525510E-05	&	13	&	8	&	1.535040E-05	&	13	&	8	 \\
 portfl-i-2    	&	87	&	1	&	13	&	0	&	1.465720E-05	&	14	&	8	&	1.460560E-05	&	14	&	8	 \\
 portfl-i-3    	&	87	&	1	&	13	&	0	&	6.280600E-06	&	14	&	8	&	6.280000E-06	&	14	&	8	 \\
 portfl-i-4    	&	87	&	1	&	13	&	0	&	2.322300E-06	&	13	&	8	&	2.297200E-06	&	13	&	8	 \\
 portfl-i-6    	&	87	&	1	&	13	&	0	&	2.524600E-06	&	13	&	8	&	2.505700E-06	&	13	&	8	 \\
 qpec1         	&	30	&	20	&	0	&	20	&	8.000000E+01	&	3	&	2	&	8.000000E+01	&	4	&	2	 \\
 qpec2         	&	30	&	20	&	0	&	20	&	4.499888E+01	&	47	&	8	&	4.499842E+01	&	206	&	8	 \\
 ralph1        	&	2	&	1	&	0	&	1	&	-2.798950E-05	&	40	&	8	&	-2.798950E-05	&	40	&	8	 \\
 ralph2        	&	2	&	1	&	0	&	1	&	-2.000000E+00	&	1	&	1	&	-2.000000E+00	&	1	&	1	 \\
 scholtes1     	&	3	&	1	&	0	&	1	&	2.000000E+00	&	6	&	2	&	2.000000E+00	&	6	&	2	 \\
 scholtes2     	&	3	&	1	&	0	&	1	&	1.500000E+01	&	5	&	2	&	1.500000E+01	&	5	&	2	 \\
 scholtes1n    	&	3	&	1	&	0	&	1	&	2.000000E+00	&	6	&	2	&	2.000000E+00	&	6	&	2	 \\
 scholtes2-n   	&	3	&	1	&	0	&	1	&	1.500000E+01	&	5	&	2	&	1.500000E+01	&	5	&	2	 \\
 scholtes3     	&	2	&	1	&	0	&	1	&	5.000000E-01	&	26	&	7	&	5.000000E-01	&	26	&	7	 \\
 scholtes4     	&	3	&	3	&	0	&	1	&	-5.590180E-05	&	46	&	8	&	-5.590180E-05	&	46	&	8	 \\
 scholtes5     	&	3	&	2	&	0	&	2	&	1.000000E+00	&	3	&	2	&	1.000000E+00	&	4	&	2	 \\
 sl1           	&	8	&	3	&	2	&	3	&	1.000001E-04	&	12	&	8	&	1.000000E-04	&	11	&	8	 \\
 stub          	&	12	&	5	&	6	&	3	&	3.283210E-05	&	28	&	2	&	3.283210E-05	&	28	&	2	 \\
 stackelberg1  	&	3	&	1	&	1	&	1	&	-3.266667E+03	&	4	&	2	&	-3.266667E+03	&	4	&	2	 \\
 tap-09        	&	86	&	36	&	32	&	32	&	1.091310E+02	&	91	&	5	&	1.091311E+02	&	136	&	8	 \\
 tap-15        	&	194	&	99	&	68	&	83	&	1.843557E+02	&	255	&	8	&	1.842949E+02	&	183	&	8	 \\
 water-net     	&	66	&	14	&	36	&	14	&	9.272644E+02	&	1317&	8   &	9.197094E+02	&	1511&	 8\\

\noalign{\smallskip} \hline
\end{tabular}
}
\end{table}

\subsection{Performance profiles}

\begin{figure}[]
\centering
\begin{minipage}[t]{0.47\linewidth}
\centering
\includegraphics[scale=0.5]{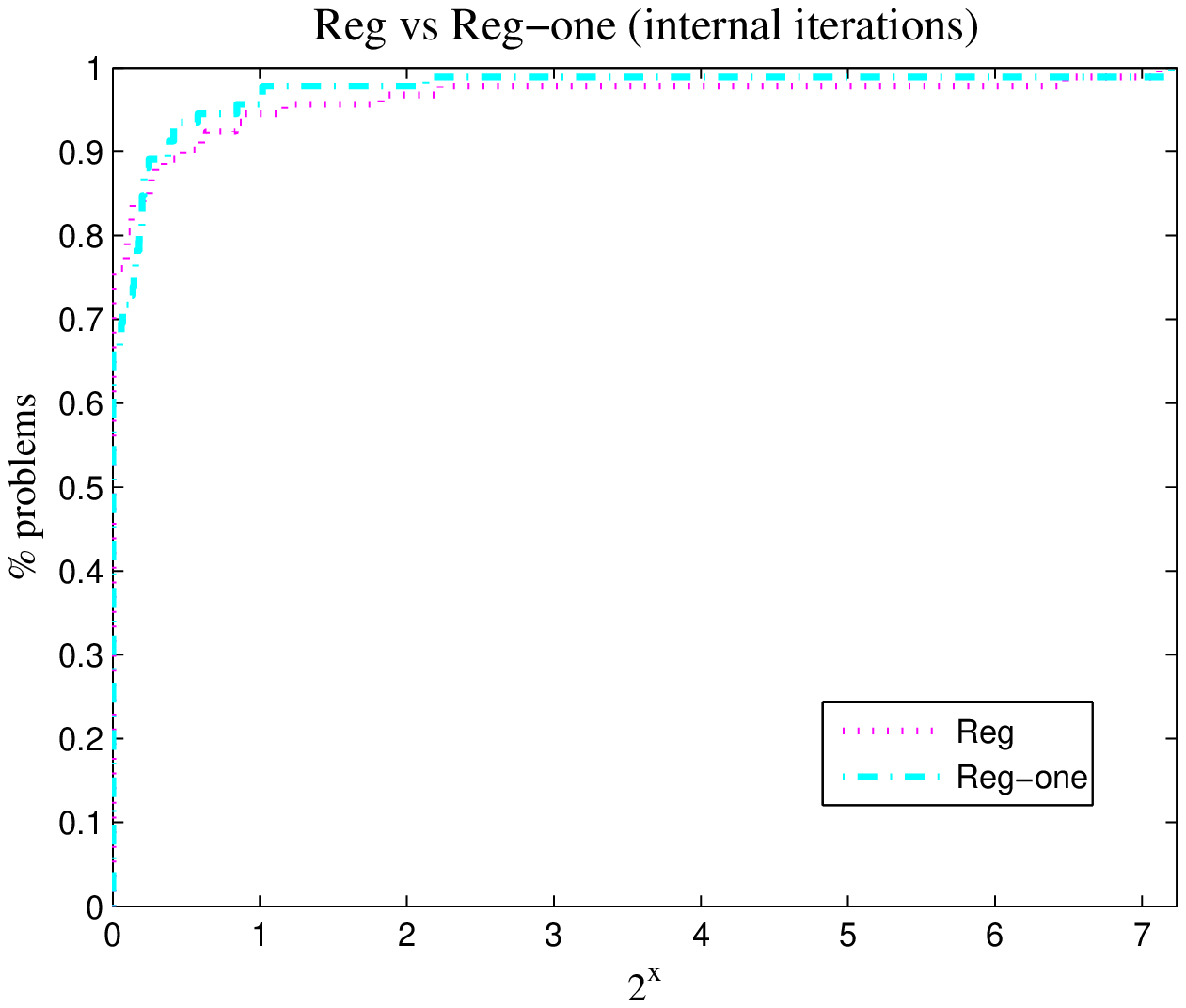}
{\caption{\label{int} \small Inner iterations performance profile}}
\end{minipage}\hspace*{\fill} \begin{minipage}[t]{0.47\linewidth}
\centering
\includegraphics[scale=0.5]{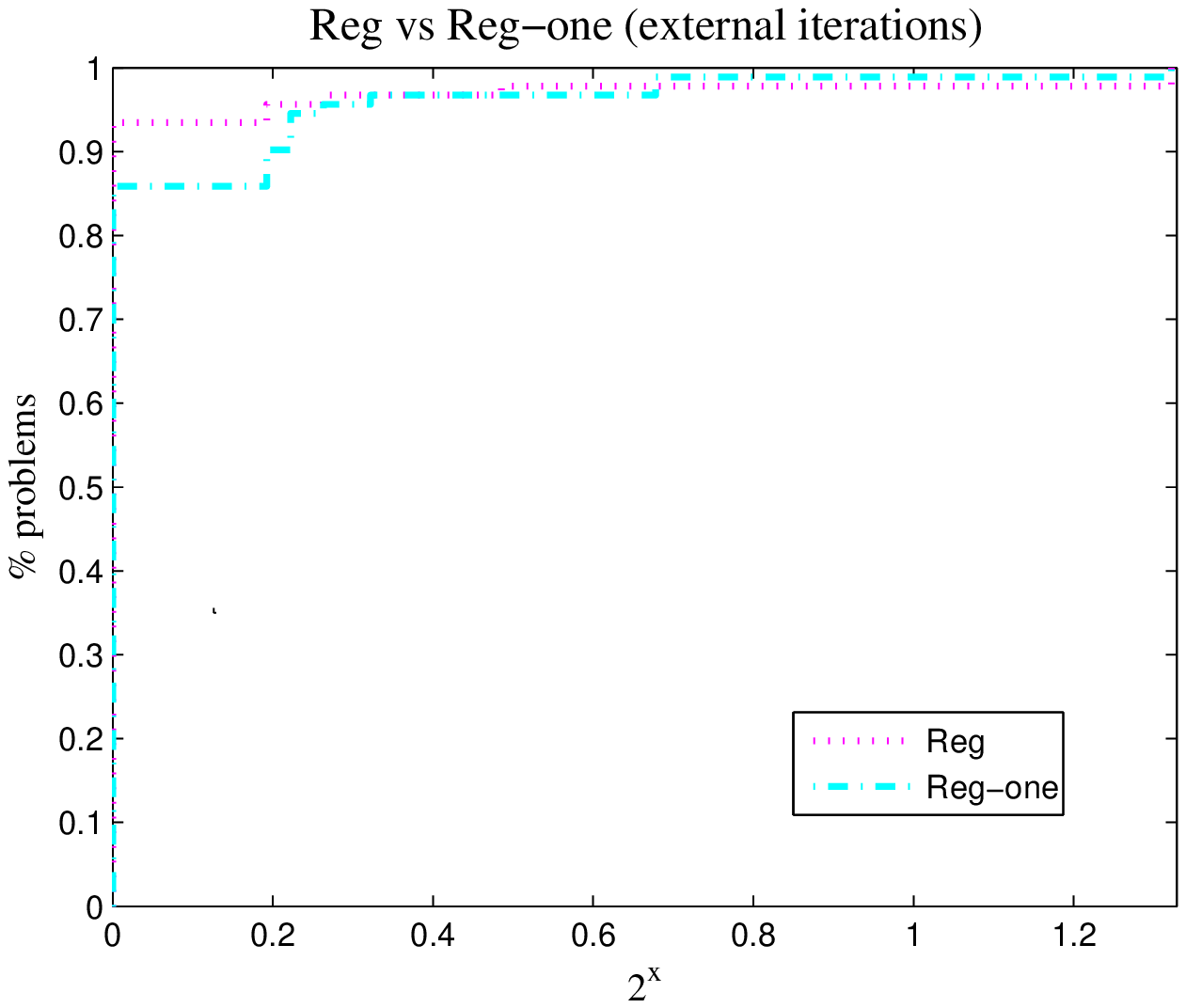}
{\caption{\label{ext} \small External iterations performance profile}}
\end{minipage}
\end{figure}

Dolan and Mor\'{e} \cite{dolan01} present a
tool to analyse the relative performance of the optimization codes
with respect to a specific metric. An easy interpretation of this graphic is that for any given metric a solver is the best
when its graphic is tending faster to 1.
The performance metrics considered are the number of
internal and external iterations, respectively. The graphics of performance profiles are in
Figure \ref{int} and Figure \ref{ext} and a $log$ scale is used.

From the Figure \ref{int} one can conclude that the Algorithm $Reg$ and $Reg$-$one$ have similar performance with respect to the inner iterations - both algorithms  present very good robustness. Concerning to the external iterations, one can attest that the Algorithm $Reg$ presents a little better performance.

\subsection{Final remarks and future work}

The first conclusion of this work is that it is possible to solve MPCC using nonlinear techniques. Four algorithms were implemented in MATLAB, combining the SQP philosophy with a regularization scheme. The numerical results of Algorithms $Reg$ and $Reg$-$one$ are very promising - both algorithms present robustness and efficiency.  The other two algorithms, $Reg$-$eq$  and $Reg$-$eq$-$one$, did not perform well. This behavior  was already expected since the equalities in the complementarity constraints still continue in both formulations. As future work a very careful analysis must be performed in order to implement some strategies to prevent this bad performance. Some important procedures like testing the algorithms with larger dimension problems and studying their convergence properties are on going.

\renewcommand{\baselinestretch}{1}
{\footnotesize
\bibliographystyle{elsart-num}

\end{document}